# The signature of a meander

Vincent Coll[1], Colton Magnant[2], Hua Wang[2]

**Abstract**

We provide a recursive classification of meander graphs, showing that each meander is identified by a unique sequence of fundamental graph-theoretic moves. This sequence is called the meander's signature. The signature not only provides a fast algorithm for the computation of the index of a Lie algebra associated with the meander, but also allows for the speedy determination of the graph's plane homotopy type - a finer invariant than the index. The signature can be used to construct arbitrarily large sets of meanders, Frobenius and otherwise, of any given size and configuration.

Making use of a more refined signature, we are able to prove an important conjecture of Gerstenhaber & Giaquinto: The spectrum of the adjoint of a principal element in a Frobenius seaweed Lie algebra consists of an unbroken chain of integers. Additionally, we show the dimensions of the associated eigenspaces to be symmetric about 0.5. In certain special cases, the signature is used to produce an explicit formula for the index of a seaweed Lie subalgebra of $\mathfrak{sl}(n)$ in terms of elementary functions.

## 1. Introduction

Meanders were introduced by Dergachev & Kirillov in [4] as planar graph representations of biparabolic (seaweed) subalgebras of $\mathfrak{sl}(n)$. In [4], the authors provided a combinatorial method of computing the index of such biparabolics from the number and type of the connected components of their associated meander graphs. Of particular interest are those seaweed algebras whose associated meander graph consists a single path. Such algebras have index zero. More generally, algebras with index zero are called *Frobenius* and have been extensively studied in the context of invariant theory [15].

[1]Department of Mathematics, Lehigh University, Bethlehem, PA, USA
[2]Department of Mathematical Sciences, Georgia Southern University, Statesboro, GA, USA.

To fix the notation, let $\mathfrak{g}$ be a Lie algebra over a field of characteristic zero. For any functional $F \in \mathfrak{g}^*$ one can associate the *Kirillov form* $B_F(x, y) = F[x, y]$ which is skew-symmetric and bilinear. The index of $\mathfrak{g}$ is defined to be the minimum dimension of the kernel of $B_F$ as $F$ ranges over $\mathfrak{g}^*$. The Lie algebra $\mathfrak{g}$ is Frobenius if its index ind $\mathfrak{g}$ is zero. Equivalently, if there is a functional $F \in \mathfrak{g}^*$ such that $B_F$ is non-degenerate. We call such an $F$ a *Frobenius functional* and the natural map $\mathfrak{g} \to \mathfrak{g}^*$ defined by $x \mapsto F([x, -])$ is an isomorphism. The image of $F$ under the inverse of this map is called a *principal element* of $\mathfrak{g}$ and will be denoted $\hat{F}$ [9]. It is the unique element of $\mathfrak{g}$ such that
$$F \circ \mathrm{ad}\hat{F} = F([\hat{F}, -]) = F.$$
When $\mathfrak{g}$ is Frobenius, we write $\mathfrak{f}$ instead of $\mathfrak{g}$. The set of Frobenius functionals of a given $\mathfrak{f}$ is, in general, quite large. The set of $F$ such that the ind $B_F$ = ind $\mathfrak{f}$ is open in $\mathfrak{f}^*$ under both the Zariski and Euclidean topologies. See [4], Proposition 2.4.

Frobenius Lie algebras are also of special interest in deformation theory and quantum group theory stemming from their connection to the classical Yang-Baxter equation (CYBE). Suppose $B_F(-,-)$ is non-degenerate and let $M$ be the matrix of $B_F(-,-)$ relative to some basis $\{x_1, \ldots, x_n\}$ of $\mathfrak{f}$. Belavin & Drinfel'd showed that $r = \sum_{i,j}(M^{-1})_{ij} x_i \wedge x_j$ is a (constant) solution of the CYBE [1]. Thus, each pair consisting of a Lie algebra $\mathfrak{f}$ together with functional $F \in \mathfrak{f}^*$ such that $B_F$ is non-degenerate provides a solution to the CYBE. See [10] and [11] for examples.

In Section 4.1, we provide a recursive classification of meander graphs, showing that each meander is identified by a unique sequence of fundamental graph-theoretic moves. This sequence is called the meander's *signature*. The signature's reduction is similar to but more sensitive than the Panyushev reduction [16] which was also used by Dvorsky [6]. The deterministic nature of the signature's *winding down* moves provides for a reverse *winding up* process which can be used to construct all meanders, Frobenius and otherwise, of any specified size and block configuration, thus providing a fertile bed of test examples. The signature not only provides a fast algorithm for the computation of the index of a Lie algebra associated with the meander, but also allows for the speedy determination of the graph's *plane homotopy type* - a finer invariant than the index.

In Section 4.2, using a *refined signature*, we show that the spectrum of the adjoint of a principal element consists of an unbroken sequence of integers $-a, -a+1, \ldots, 0, 1, \ldots, a, a+1$. This proves an outstanding conjecture [11] of Gerstenhaber & Giaquinto. In fact, this was stated as a result in [9] and [11] but it was later observed by Duflo [5], in a personal communication to Gerstenhaber & Giaquinto, to contain a gap in the proof. Additionally, we show the dimensions of the associated eigenvalues to be symmetric about 0.5. In a personal communication to the first author, Giaquinto et al. [12] offer that a different approach using an alternative Frobenus functional (a generalized cyclic functional [11]) may also be used to show that the spectrum is unbroken and symmetric.

In Section 4.3, we use the notion of the signature to prove a formula for the index of a meander with 4 blocks. This is a natural extension of results of Elashvili [7] and Coll et al. [2]. Section 4.4 contains another application of the signature whereby we produce new infinite families of (parabolic and bi-parabolic) Frobenius meanders with arbitratily large blocks and of arbitrarily large size. Finally in Section 5, we present some conjectures stemming from this work.

## 2. Meanders

We begin with the definition of a seaweed subalgebra, then we connect this to the Frobenius functional of Dergachev & Kirillov in [4] through the underlying graph of the directed graph of the functional.

**Definition 1.** *Let $\mathbf{k}$ be an arbitrary field of characteristic $0$ and $n$ a positive integer. Fix two ordered partitions $\{a_i\}_{i=1}^k$ and $\{b_j\}_{j=1}^\ell$ of the number $n$. Let $\{e_i\}_{i=1}^n$ be the standard basis in $\mathbf{k}^n$. A subalgebra of $\mathrm{M}_{\mathrm{n} \times \mathrm{n}}(\mathbf{k})$ that preserves the vector spaces $\{V_i = \mathrm{span}(\mathrm{e}_1, \ldots, \mathrm{e}_{\mathrm{a}_1 + \cdots + \mathrm{a}_\mathrm{i}})\}$ and $\{W_j = \mathrm{span}(\mathrm{e}_{\mathrm{b}_1 + \cdots + \mathrm{b}_\mathrm{j} + 1}, \ldots, \mathrm{e}_\mathrm{n})\}$ is called a subalgebra of* seaweed *type due to the suggestive shape of the subalgebra in the total matrix algebra. See Figure 1 where the left matrix algebra is a seaweed subalgebra of $\mathfrak{sl}(3)$ with $a_1 = 1$, $a_2 = 2$ and $b_1 = 3$. The right matrix algebra is a seaweed subalgebra of $\mathfrak{sl}(7)$ with $a_1 = a_2 = 2$, $a_3 = 3$, $b_1 = 5$ and $b_2 = 2$.*

To each seaweed algebra, a planar graph, called a *meander*, can be constructed as follows: Line up $n$ vertices and label them with $v_1, v_2, \ldots, v_n$ and partition this set into two ordered partitions, called *top* and *bottom*. For each part in the top (likewise bottom) we build up the graph by adding edges in the same way. This involves adding the edge from the first vertex of a part to the last vertex of the same part drawn concave down (respectively concave up in the bottom part case). The edge addition is then repeated between the second vertex and the second to last and so on within each part of both partitions. Note that this is simply the underlying (undirected) graph of the directed graph of the functional. The parts of these partitions are called *blocks*. With top blocks of sizes $a_1, a_2, \ldots, a_\ell$ and bottom blocks of sizes $b_1, b_2, \ldots, b_m$, we say such a meander has *type* $\frac{a_1|a_2|\ldots|a_\ell}{b_1|b_2|\ldots|b_m}$. See Figure 2 where the left meander is of type $\frac{1|2}{3}$ and the right meander is of type $\frac{2|2|3}{5|2}$.

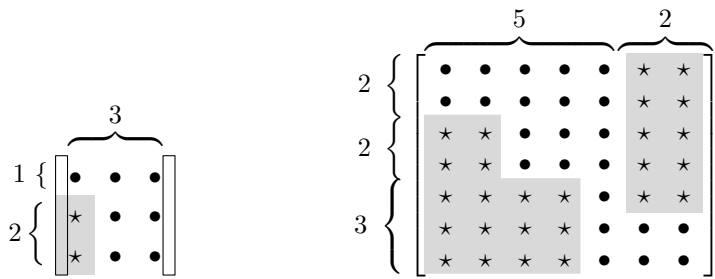

Figure 1: Seaweed algebras.

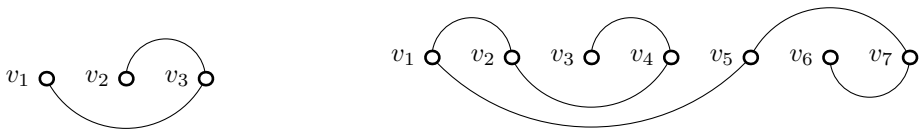

Figure 2: Meanders.

To each Frobenius seaweed subalgebra, there is a canonically defined functional developed by Dergachev & Kirillov in [4]. To present the definition of the functional, we require some preliminary notation. Let $e_{ij}$ denote the $n \times n$ matrix with a 1 in the $(i,j)$ location and 0's everywhere else. If $S$ is a subset of the indices $(i,j)$ with $1 \leq i, j \leq n$, then $F_S$ denotes the functional $\sum_{s \in S} e_s^*$ which acts on the space $M_n$ of all $n \times n$ matrices. We say that $S$ *supports* $F_S$ and the *directed graph of the functional*, $\Gamma(S)$, has vertices $v_1, \ldots, v_n$ with an arrow from $v_i$ to $v_j$ whenever $(i,j) \in S$. See Figures 3.

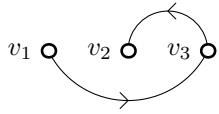

Figure 3: $\Gamma(S)$ where $S = \{(1,3),(3,2)\}$ and $F_S = e_{13}^* + e_{32}^*$.

Note that the underlying (undirected) graph of $\Gamma(S)$ is the meander of type $\frac{1|2}{3}$. See the left meander in Figure 2. In this case, we have

$$F = \begin{bmatrix} 1 & 2 & 1 \\ 0 & -1 & -1 \\ 0 & 1 & 0 \end{bmatrix} \quad \text{and} \quad \hat{F} = \begin{bmatrix} 1 & 0 & 0 \\ 0 & -1 & 0 \\ 0 & 0 & 0 \end{bmatrix}.$$

## 3. Recursive classification and winding down

In this section, we show that any meander can be contracted ("Wound Down") to the empty meander through a sequence of graph-theoretic moves; each of which is uniquely determined by the structure of the meander at the time of move application. We find that there are five such moves, only one of which affects the component structure of the graph and is therefore the only move capable of modifying the *index* of the graph, here defined to be twice the number of cycles plus the number of paths minus one. Dergachev & Kirillov [4] showed that the index of the graph is precisely the index of the associated seaweed subalgebra. Since the sequence of moves which contracts a meander to the empty meander uniquely identifies the graph, we call this sequence the meander's *signature*.

**Remark** Although this reduction is very similar to the reduction of Panyushev [16] and Dvorsky [6], we define the operations carefully so as to facilitate reversing the process. This allows us to construct all meanders of any size, number of blocks or index by a simple sequence of moves.

We recall the definition of a contraction from basic graph theory. The *contraction* of an edge $e = uv$ is the removal of $e$ and identification of $u$ and $v$ to form a single vertex. By "$\equiv$", we mean that two meanders have the same component structure which we later call the meander graph's plane homotopy type.

The move performed in Lemma 2 is determined based on information about the relationship between the sizes of $A_1$ and $B_1$, the first top and first bottom blocks respectively. In most cases, edges of $A_1$ are contracted, sometimes leaving other vertices untouched, sometimes rotating other vertices in the process. At times, these contractions are among the outer edges of $A_1$ but at other times, the case requires that all edges of $A_1$ be contracted.

We state Lemma 2 without proof because we later prove a more refined version of this lemma (see Lemma 5) for use in the proof of our main result concerning the spectrum.

**Lemma 2 (Winding Down - Simplified).** *Let $M$ be a general meander of type*

$$M = \frac{a_1|a_2|a_3|\ldots|a_n}{b_1|b_2|b_3|\ldots|b_m}$$

*Then we get the following cases:*

1. **Flip ($F_0$)** *If $a_1 < b_1$, then simply exchange $a_i$ for $b_i$ to get*

   $$M \equiv \frac{b_1|b_2|\ldots|b_m}{a_1|a_2|\ldots|a_n}.$$

2. **Component Elimination ($C_0(c)$)** *If $a_1 = b_1$, then*

   $$M \to \frac{a_2|a_3|\ldots|a_n}{b_2|b_3|\ldots|b_m}.$$

3. **Block Elimination ($B_0$)** *If $a_1 = 2b_1$, then*

   $$M \equiv \frac{b_1|a_2|a_3|\ldots|a_n}{b_2|b_3|\ldots|b_m}.$$

4. **Rotation Contraction ($R_0$)** *If $b_1 < a_1 < 2b_1$, then*

   $$M \equiv \frac{b_1|a_2|a_3|\ldots|a_n}{(2b_1 - a_1)|b_2|b_3|\ldots|b_m}.$$

5. **Pure Contraction ($P_0$)** *If $a_1 > 2b_1$, then*

   $$M \equiv \frac{(a_1 - 2b_1)|b_1|a_2|a_3|\ldots|a_n}{b_2|b_3|\ldots|b_m}.$$

We call these cases *moves* since they move one meander to another. Notice that in each of these moves except the Flip, the number of vertices in the graph is reduced. Also note that in each move except the Component Elimination, the component structure of the meander is maintained.

Given a meander, there exists a unique sequence of moves (elements of the set $\{F_0, C_0(c), B_0, R_0, P_0\}$) to reduce down to the empty meander. Such a list is called the *signature* of the meander. Given a different set of moves, as in Lemma 5, there will naturally be a different signature.

As examples of these moves, we provide the following figures. Let $M'$ be the resulting meander after the move is performed on $M$. Figure 4 demonstrates a Block Elimination move where the dotted edges are in $M$ but not $M'$ while the bold edges are in $M'$ but not $M$. Here $a_1 = 10$ and $b_1 = 5$. Edges and vertices of the other blocks are not drawn (here we are focusing on the blocks $A_1$ and $B_1$). Figure 5 presents a Rotation Contraction where $a_1 = 10$, $b_1 = 7$ and all other blocks are not drawn. Figure 6 demonstrates a Pure Contraction where $a_1 = 11$ and $b_1 = 4$ and all other blocks are not drawn.

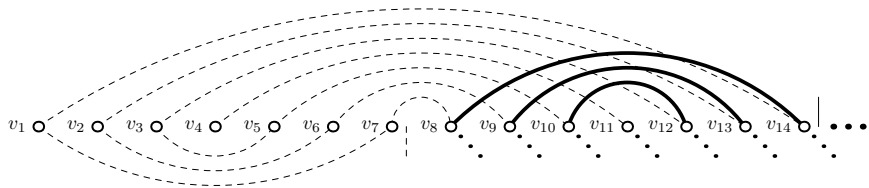

Figure 4: Block Elimination.

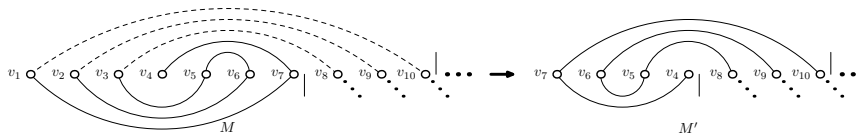

Figure 5: Rotation Contraction.

Lemma 2 provides a fast algorithm for computing the signature of a meander. Indeed, from the signature, we can read off the index of the meander by simply adding the parameters used in the Component Elimination moves and subtracting 1. That is, if $\{c_1, c_2, \ldots, c_m\}$ are the parameters used in the Component Elimination moves, then the index of the meander equals $\sum_{i=1}^{m}(c_i) - 1$.

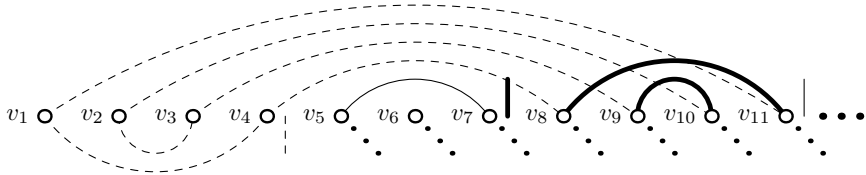

Figure 6: Pure Contraction.

**Corollary 3.** *There is a linear (in the order of the meander) time algorithm for computing the index of the meander.*

There is more information in the signature of a meander than simply the index. For example, a meander whose signature contains $C_0(3)$ is fundamentally different from a meander whose signature contains $C_0(1)$ and $C_0(2)$ separately. Geometrically, these situations can be reduced down to a compressed form, the picture of the component at the instant we perform the Component Elimination move. In the case of the $C_0(3)$, this is a circle with a point in the middle while the case of $C_0(1)$ and $C_0(2)$ is a circle with a point outside. See Figure 7. We call these geometric figures the *symbols* of the Component Elimination moves and the collection of symbols is called the graph's *plane homotopy type*. This provides a more refined invariant than the index of the meander. Indeed, the meanders in Figure 7 both have the same index, namely 2.

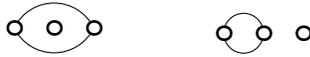

Figure 7: Meanders $\frac{3}{3}$ and $\frac{2|1}{2|1}$.

In particular, Figure 8 demonstrates the plane homotopy type of a meander whose signature includes $C_0(7)$, $C_0(5)$, $C_0(4)$, $C_0(4)$, $C_0(3)$, $C_0(2)$, $C_0(1)$, $C_0(1)$. Note that the order of these symbols is not important and any of these could be an Internal Component Elimination (defined in Lemma 5) as opposed to the Component Elimination.

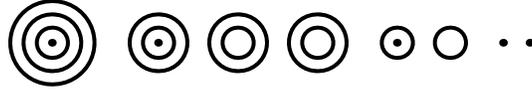

Figure 8: Plane homotopy type.

Since the only move that changes the component structure is $C_0$, we get the following corollary which provides a recursive classification of Frobenius meanders.

**Theorem 4 (Recursive Classification).** *A meander is Frobenius if and only if the signature contains no $C_0$ except the very last move which must be $C_0(1)$.*

As an example of a Frobenius meander, consider the meander of type $\frac{6|1}{2|3|2}$. This meander has signature $P_0 F_0 R_0 B_0 F_0 B_0 F_0 B_0 C_0(1)$ which is demonstrated in Figure 9.

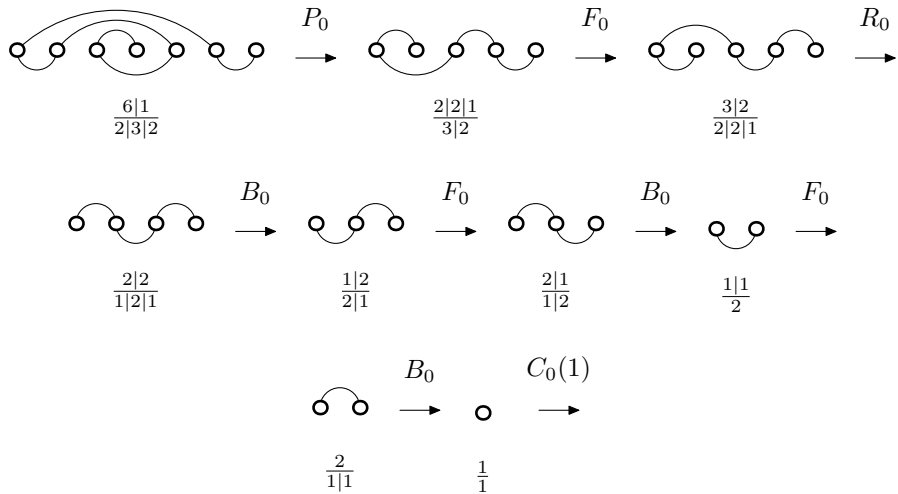

Figure 9: Detailed signature of $\frac{6|1}{2|3|2}$.

As an example of a non-Frobenius meander, consider the meander of type $\frac{16|2|4}{5|17}$. This meander has signature $P_0 F_0 P_0 C_0(5) P_0 F_0 B_0 C_0(2)$ which is demonstrated in Figure 10.

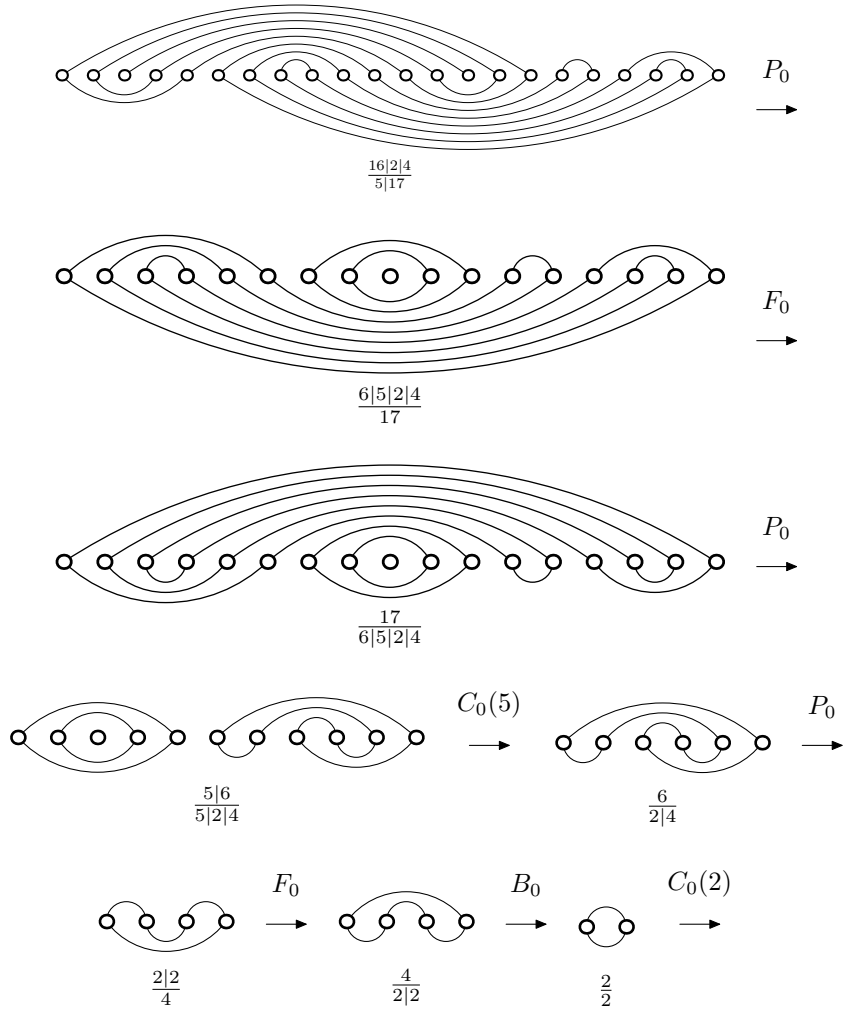

Figure 10: Detailed signature of $\frac{16|2|4}{5|17}$.

## 4. Refined signature and applications

As mentioned in Section 3, there are certainly simple sets of Winding Down moves that would suffice to classify a meander. We use this set of moves for the purpose of proving important results about the dimensions of eigenvalues in Subsection 4.2. Note that the Pure Contraction from Lemma 2 has been changed and we have added Internal moves.

**Lemma 5 (Winding Down).** *Let $M$ be a general meander of the form*

$$M = \frac{a_1|a_2|a_3|\ldots|a_n}{b_1|b_2|b_3|\ldots|b_m}.$$

*Then we get the following cases:*

1. **Flip ($F$)** *If $a_1 < b_1$, then flip the meander vertically to obtain*

$$M \equiv \frac{b_1|b_2|\ldots|b_m}{a_1|a_2|\ldots|a_n}.$$

2. **Component Elimination ($C(c)$)** *If $a_1 = b_1 = c$ for some number $c$, then*

$$M \to \frac{a_2|a_3|\ldots|a_n}{b_2|b_3|\ldots|b_m}.$$

3. **Block Elimination ($B$)** *If $a_1 = 2b_1$, then*

$$M \equiv \frac{b_1|a_2|a_3|\ldots|a_n}{b_2|b_3|\ldots|b_m}.$$

4. **Rotation Contraction ($R$)** *If $b_1 < a_1 < 2b_1$, then*

$$M \equiv \frac{b_1|a_2|a_3|\ldots|a_n}{(2b_1 - a_1)|b_2|b_3|\ldots|b_m}.$$

5. **Internal Component Elimination ($IC(c)$)** *Suppose $a_1 > 2b_1$ and the center point of $A_1$ is also the center of a bottom block $B_i$. Then $c = b_i$ and*

$$M \to \frac{(a_1 - b_i)|a_2|a_3|\ldots|a_n}{b_1|b_2|\ldots|b_{i-1}|b_{i+1}|b_{i+2}|\ldots|b_m}.$$

6. **Internal Block Elimination ($IB$)** *Suppose $a_1 > 2b_1$, $a_1$ is even and the vertices $v_{a_1/2}$ and $v_{a_1/2+1}$ are end vertices of bottom blocks $B_i$ and $B_{i+1}$ respectively for some $i$. Then*

$$M \equiv \frac{(a_1 - b_i)|a_2|a_3|\ldots|a_n}{b_1|b_2|\ldots|b_{i-1}|b_{i+1}|b_{i+2}|\ldots|b_m}.$$

7. **Internal Rotation Contraction ($IR$)** *Suppose $a_1 > 2b_1$ and either $a_1$ is odd or the vertices $v_{a_1/2}$ and $v_{a_1/2+1}$ are in the same bottom block $B_i$ for some $i$ where the center of $A_1$ is not also the center of $B_i$. Let* $t = \lfloor a_1/2 \rfloor$ and let $B_i$ be the block containing $v_t$. Let $r$ be the shortest distance between an end of $B_i$ and the center of $A_1$ (note that $r$ need not be an integer). Let $s = 2r + 1$. Then

$$M \equiv \frac{(a_1 - b_i + s)|a_2|a_3|\ldots|a_n}{b_1|b_2|\ldots|b_{i-1}|s|b_{i+1}|b_{i+2}|\ldots|b_m}.$$

It can be easily verified that these cases are disjoint and cover all the possibilities.

*Proof:* Let $v_1 v_2 \ldots v_{a_1}$ be the vertices of the first top block in $M$. We consider each case separately.

**Case 1 (Flip).** $a_1 < b_1$.

This case is trivial since we merely move all top edges to the bottom and all bottom edges to the top.

**Case 2 (Component Elimination).** $a_1 = b_1$.

In this case, we create a new meander $M'$ by removing the top block $A_1$ and the corresponding bottom block $B_1$ to obtain $M' = \frac{a_2|a_3|\ldots|a_n}{b_2|b_3|\ldots|b_m}$. The Component Elimination cases (and Internal Component Elimination) are the only cases in which the component structure of the meander is changed since components are removed.

**Case 3 (Block Elimination).** $a_1 = 2b_1$.

Create a new meander $M'$ by contracting all edges of the form $v_i v_{a_1-i+1}$ where $1 \leq i \leq b_1$. Note that, since $a_1 = 2b_1$, all top edges in the first set have been contracted by this process. Here, the first set of bottom edges in $M$ becomes the first set of top edges in $M'$. See Figure 4. Clearly $M \equiv M'$, as desired.

**Case 4 (Rotation Contraction).** $b_1 < a_1 < 2b_1$.

Create a new meander $M'$ by contracting all edges of the form $v_i v_{a_1-i+1}$ where $1 \leq i \leq a_1 - b_1$. Also exchange the positions of all pairs of vertices of the form $v_i v_{a_1-i+1}$ where $a_1 - b_1 < i < \frac{a_1}{2}$. Here, the $b_1$ vertices of the first bottom block are rotated clockwise around the point $\frac{a_1+1}{2}$. If $a_1$ is odd, this point is a vertex and that vertex remains where it was. If $a_1$ is even, then the point is between two vertices and is simply a center of rotation.

In this case, some of the top edges in the first set of $M$ get contracted while the remaining top edges in this first set get switched to bottom edges. This leaves $(2b_1 - a_1)$ vertices in the first bottom set in $M'$. The bottom edges of the first set of $M$ get switched to top edges and so we have $b_1$ vertices in the first top set in $M'$. See Figure 5.

In Figure 5, the vertices $v_1, v_2$ and $v_3$ have been identified with $v_{10}, v_9$ and $v_8$ respectively. Also the center vertices $v_4, v_5, v_6$ and $v_7$ have rotated $180°$ around the point between $v_5$ and $v_6$ to reverse their order. The three edges that were on the bottom have rotated to the top while the two edges that were the innermost two on the top have rotated to the bottom. Thus, $M \equiv M'$.

**Case 5 (Internal Component Elimination).** $a_1 > 2b_1$ and the center point of $A_1$ is also the center of a bottom block $B_i$

In this case, we create a new meander $M'$ by removing the block $B_i$ and the corresponding edges of $A_1$ to obtain $M' = \frac{(a_1 - b_i)|a_2|a_3|...|a_n}{b_1|b_2|...|b_{i-1}|b_{i+1}|b_{i+2}|...|b_m}$. Note that the component structure of the meander is changed since components are removed.

**Case 6 (Internal Block Elimination).** $a_1 > 2b_1$, $a_1$ is even and the vertices $v_{a_1/2}$ and $v_{a_1/2+1}$ are end vertices of bottom blocks $B_i$ and $B_{i+1}$ respectively for some $i$.

Since this case is similar to the Block Elimination case, Case 3 above, we present only a diagram of this case. In Figure 11, $a_1 = 16$, $b_1 = 3$, $b_2 = 5$ and $b_3 = 4 < 8 = a_1/2$. The rest of the blocks are not pictured.

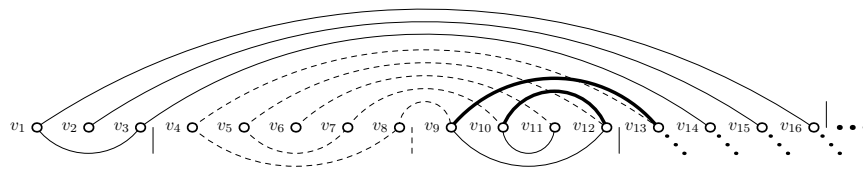

Figure 11: Reduction in Case 6. Internal Block Elimination.

**Case 7 (Internal Rotation Contraction).** $a_1 > 2b_1$ *and either $a_1$ is odd or the vertices $v_{a_1/2}$ and $v_{a_1/2+1}$ are in the same bottom block $B_i$ for some $i$ where the center of $A_1$ is not also the center of $B_i$.*

This case is similar to the Rotation Contraction case, Case 4 above. In Figure 12, $a_1 = 14$, $b_1 = 3$, and $b_2 = 6$. The rest of the blocks are not pictured.

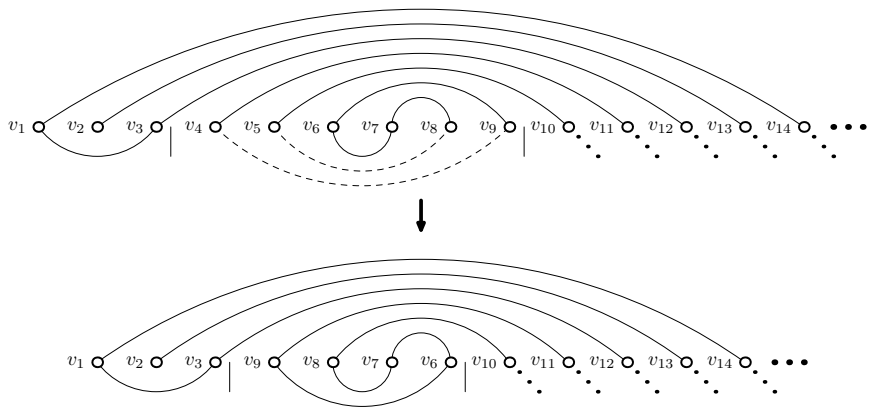

Figure 12: Reduction in Case 7. Internal Rotation Contraction.

Intuitively, let $g$ be the gap between bottom blocks that is closest to the center of the $A_1$ block (between $v_9$ and $v_{10}$ in Figure 12). Going from $g$, if we follow the top edges to the gap between the two vertices on the opposite side of the center of $A_1$ (between $v_5$ and $v_6$ in Figure 12), then the bottom edges below that gap will be contracted.

More formally, suppose $a_1$ is odd and let $m = (a_1 + 1)/2$. Let $B$ be the bottom block containing $v_m$ where $B$ contains vertices from $v_j$ to $v_k$ (inclusive) and suppose that $(m - j) > (k - m)$. All other cases, including the case where $a_1$ is even, proceed similarly. The move involves contracting the bottom edges incident to vertices $v_x$ where $j \leq x \leq 2m - k - 1$.

This completes the proof of this case and the proof of Theorem 5. □

*4.1. Winding up*

By Lemma 5, the following result is immediate by simply reversing the contraction process.

**Lemma 6 (Winding Up).** *Every meander is the result of a unique sequence of the following moves starting from the empty meander. Starting with a meander $M_i = \frac{a_1|a_2|...|a_n}{b_1|b_2|...|b_m}$, create a meander $M_{i+1}$ by one of of the following:*

1. **Flip ($\tilde{F}$)** *Exchange the list of $a_j$ for the list of $b_j$ (namely flip the entire meander vertically).*
2. **Component Creation with parameter $c$ ($\tilde{C}(c)$)** *To the left of the block $B_1$ on the bottom and to the left of the block $A_1$, add a new block of size $c$.*
3. **Block Creation ($\tilde{B}$)** *To the left of the block $B_1$ on the bottom, add a new block of size $a_1$ and replace the block $A_1$ with a block of size $2a_1$. This creates $M_{i+1}$ of type $\frac{2a_1|a_2|...|a_n}{a_1|b_1|b_2|...|b_m}$.*
4. **Rotation Expansion ($\tilde{R}$)** *Replace the block $B_1$ on the bottom with a block of size $a_1$ and replace the block $A_1$ on the top with a block of size $b_1 + 2(a_1 - b_1)$. This creates $M_{i+1}$ of type $\frac{(2a_1-b_1)|a_2|...|a_n}{a_1|b_2|...|b_m}$. Note that this expansion requires that $a_1 > b_1$.*
5. **Internal Component Creation with parameter $c$ ($I\tilde{C}(c)$)** *This move can only be made if $a_1$ is even and the vertex $v_{a_1/2}$ is the right-most vertex of some bottom block $B_j$. Then insert a block $C$ of order $c$ between $B_j$ and $B_{j+1}$, adding $c$ vertices to $A_1$ to create the meander $M_{i+1}$ of type $\frac{(a_1+c)|a_2|a_3|...|a_n}{b_1|b_2|...|b_j|c|b_{j+1}|b_{j+2}|...|b_m}$.*
6. **Internal Block Creation with parameter $b \geq 2$ ($I\tilde{B}(b)$)** *This move can only be made if $a_1 > 2b_1$ and $B_b$ is a bottom block with left-most vertex $v_k$ where $k \leq a_1/2$. Between $B_{b-1}$ and $B_b$, create a bottom block of order $a_1 - 2(k-1)$ and increase $a_1$ accordingly to create the meander $M_{i+1}$ of type $\frac{(2a_1-2(k-1))|a_2|a_3|...|a_n}{b_1|b_2|...|b_{b-1}|(a_1-2(k-1))|b_b|b_{b+1}|...|b_m}$.*
7. **Internal Rotation Expansion ($I\tilde{R}$)** *This move can only be made if $a_1 > 2b_1$, the center of $A_1$ is not also the center of a bottom block and if $a_1$ is even, the point $v_{a_1/2}$ is not the right-most vertex of a bottom block. Let $B_j$ be the bottom block containing the center point of $A_1$ and let $v_k$ and $v_\ell$ be the end-vertices of $B_j$ that are farthest and closest respectively to the center of $A_1$. Then set $r = |k - a_1/2| - |a_1/2 - \ell|$ and create the meander $M_{i+1}$ of type $\frac{(a_1+r)|a_2|a_3|...|a_n}{b_1|b_2|...|b_{j-1}|(b_j+r)|b_{j+1}|b_{j+2}|...|b_m}$.*

Starting with the empty meander, the first operation must be the Component Creation move and, if the meander is to be Frobenius, this must be the only Component Creation move used and it has parameter $c = 1$. Thus, the following exciting result is immediate.

**Theorem 7.** *Every Frobenius meander can be constructed from the empty meander by a sequence of Winding Up moves where the first move must be $\tilde{C}(1)$ and the moves $\tilde{C}$ or $I\tilde{C}$ never appear again in the sequence.*

In fact, any meander (Frobenius or otherwise) can be constructed by a sequence of Winding Up moves.

**Theorem 8.** *Every meander can be constructed from the empty meander by a sequence of Winding Up moves (where the first move must be $\tilde{C}(c)$ for some $c$).*

Using the simplified signature of Lemma 2, we may also define a simplified "Winding Up" sequence. Such a sequence consists of moves labeled $\tilde{F}_0, \tilde{C}_0(c), \tilde{B}_0, \tilde{R}_0$ and $\tilde{P}_0$. For example, from the sequence $\tilde{C}_0(2)\tilde{B}_0\tilde{F}_0\tilde{P}_0\tilde{C}_0(5)\tilde{P}_0\tilde{F}_0\tilde{P}_0$, we get the meander $\frac{16|2|4}{5|17}$ as seen in Figure 10 by reversing the arrows and putting hats on all the move labels.

*4.2. Spectrum*

Let $\hat{F}$ be a principle element of a Frobenius seaweed Lie algebra. In this section, we compute the spectrum of $\hat{F}$. To do this, we render the eigenvalue algorithm of [11] into a graph-theoretic definition of eigenvalues of a meander.

For our main result concerning the eigenvalues of $\text{ad}\hat{F}$, we render the eigenvalue algorithm of [11] into a graph-theoretic definition of eigenvalues of a meander. This is defined as follows. Let $M$ be a meander drawn in the standard way and orient all top edges to the left and all bottom edges to the right. Define the *measure* of a path from $u$ to $v$ in this meander to be the number of forward edges minus the number of backward edges encountered when moving from $u$ to $v$ on this path and we also call this the measure of the pair $(u, v)$. We call a pair of vertices $(v_i, v_j)$ *admissible* if one of the following cases holds: $i < j$ and $v_i$ and $v_j$ are in the same bottom block or $i \geq j$ and $v_i$ and $v_j$ are in the same top block. Over all admissible pairs $(v_i, v_j)$, calculate the measure of the $(v_i, v_j)$ path through the meander. This measure is an eigenvalue. The set of all such measures is precisely the set of eigenvalues where the multiplicity of the measures in the list is the dimension of the eigenvalue where we subtract one from the dimension of 0. Note that,

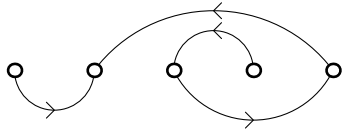

Figure 13: The oriented meander $\frac{1|4}{2|3}$.

from this algorithm, it is clear that this spectrum contains only integers as mentioned in [11].

As an example of the computation of this spectrum, consider the meander $\frac{1|4}{2|3}$. With the orientations, this is the meander pictured in Figure 13. To represent the admissible pairs versus the non-admissible pairs, we make a $5 \times 5$ matrix where the bold elements will represent the admissible pairs.

$$\begin{bmatrix} \mathbf{0} & \mathbf{1} & 0 & 0 & 0 \\ 0 & \mathbf{0} & 0 & 0 & 0 \\ 0 & \mathbf{2} & \mathbf{0} & -1 & \mathbf{1} \\ 0 & \mathbf{3} & \mathbf{1} & \mathbf{0} & \mathbf{2} \\ 0 & \mathbf{1} & -\mathbf{1} & -\mathbf{2} & \mathbf{0} \end{bmatrix}$$

This matrix is simply a bookkeeping device for recording the dimensions of the eigenvalues but it also has the same block diagonal structure as the seaweed subalgebra itself. This shows that the spectrum is $\{-2, -1, 0, 1, 2, 3\}$ with corresponding dimensions $\{1, 2, 4, 4, 2, 1\}$ respectively.

For ease of notation, call the eigenvalues of a principle element associated with the meander the eigenvalues of the meander. We call a spectrum *symmetric* if it is symmetric about 0.5, namely ranging from $-a$ to $a+1$ for some integer $a$ where the dimension of $e$ equals the dimension of $-e+1$ for all $-a \leq e \leq a+1$. We call a spectrum *unbroken* if the eigenvalues form an unbroken sequence of integers.

The following theorem confirms a conjecture of Gerstenhaber & Giaquinto [11] under the additional restriction that the dimensions of the eigenvalues are symmetric about 0.5.

**Theorem 9.** *The spectrum of a Frobenius meander is symmetric and unbroken.*

*Proof:* For a set of vertices $A$ (where $A$ will usually be a block), let $\sigma(\overleftarrow{A})$ denote the list of measures with multiplicity contributed by admissible pairs in $A$ where the first element of each pair is to the right of the second element. To this set, we add the element 0 with multiplicity $\left\lfloor \frac{|A|}{2} \right\rfloor$ to count the pairs of a vertex with itself (note that each vertex is counted once in the a block and once in a bottom block). Similarly define $\sigma(\overrightarrow{A})$. Given $\sigma(A)$ (with some direction), define $\sigma(A)+1$ to be the list of measures where 1 is added to each eigenvalue. Along the same lines, define $-\sigma(A)$ to be the list of negatives of the measures in $\sigma(A)$.

In order to prove this statement, we use the eigenvalue algorithm above along with the Winding Up Lemma, Lemma 6. The proof is by induction on the number of moves performed in the Winding Up process. The spectrum of the meander consisting of a single vertex is empty so it is trivially symmetric and unbroken so we may assume $M$ is a Frobenius meander for which the spectrum is symmetric and unbroken on at least two vertices.

Suppose also that the measures contributed within each block of the meander (in the appropriate directions) are symmetric and unbroken with center 0.5. For example, the set of measures contributed by the block $A_2$ (the block of size 4) in Figure 13 contributes $(-2, -1, 0, 0, 1, 1, 2, 3)$. In proving this theorem, we are thereby proving the following fact.

**Fact 10.** *In any Frobenius meander, the measures contributed within each block (in the appropriate direction) are symmetric and unbroken with center 0.5.*

Performing a Flip move on this meander has no effect on the spectrum. Also the Component Creation and Internal Component Creation will never be used since we consider only Frobenius meanders. Thus, we break the proof into five cases based on which of the remaining moves is performed to produce a larger meander $M'$ from $M$.

The proof of each case proceeds in the same general way. We consider cases based on which Winding Up operation is performed to get from a meander $M$ to a new larger meander $M'$. In each case, we use the measures contributed by groups of vertices in $M$ to construct the measures contributed by groups in $M'$. The measures contributed by blocks not involved in the moves are clearly not affected by the move so these measures are not mentioned in the argument. Below, we prove the first case carefully and, since all other cases follow similarly, we simply list the subcases for each of the other cases. One of the induction assumptions, namely the one concerned with contributions of measures between pairs of blocks, does not pose a challenge until Case 3. We therefore provide some insight into the proof of Case 3.

Within the cases below, the subcases are simply listing the sets of measures that have changed in the Winding Up process. We list them only for the purpose of referring to those that pair off, thereby making the resulting full set of measures symmetric and unbroken.

**Case 1.** *A Block Creation move is performed.*

Let $A$ be the set consisting of the first $a_1$ vertices of $M$. Let $A'$ be the set consisting of the first $a_1$ vertices of $M'$ and let $B'$ be the set consisting of the next $a_1$ vertices of $M'$, namely $v'_{a_1+1}, \ldots, v'_{2a_1}$. A simple consideration of each pair within these segments yields the following.

1. $\sigma(\overrightarrow{A'}) = \sigma(\overleftarrow{A})$.
2. $\sigma(\overleftarrow{A'}) = -\sigma(\overleftarrow{A})$.
3. $\sigma(\overleftarrow{B'}) = \sigma(\overleftarrow{A})$.
4. $\sigma(B' \to A') = [\sigma(\overleftarrow{A}) + 1] \cup [-\sigma(\overleftarrow{A}) + 1]$.

The measures contributed in Items 2 and 3 are balanced by those contributed in the first and second parts respectively of Item 4. Since the contribution of each block in $M$ was symmetric and unbroken and the pairs of segments that balance each other are in the same block, the contribution of each block to the measures in $M'$ is still symmetric and unbroken as desired.

**Case 2.** *A Rotation Expansion move is performed.*

In order for a Rotation Expansion move to be performed, we must have $a_1 > b_1$. Let $A$ be the first $b_1$ vertices of $M$, namely $v_1, v_2, \ldots, v_{b_1}$ and let $B$ be what remains of the first $a_1$ vertices of $M$, namely $v_{b_1+1}, \ldots, v_{a_1}$. After the Rotation Expansion is performed, let $A'$ be the first $a_1 - b_1$ vertices of $M'$, let $B'$ be what remains of the first $a_1$ vertices of $M'$ and let $C'$ be the next $a_1 - b_1$ vertices of $M'$, namely $v'_{a_1+1}, \ldots, v'_{a_1+b_1}$.

A simple consideration of each pair within these segments yields the following, in which the sigma function of each set in $M'$ is defined in terms of the sigma function of a corresponding set in $M$.

1. $\sigma(\overrightarrow{A'}) = \sigma(\overleftarrow{B})$.

2. $\sigma(\overleftarrow{A'}) = -\sigma(\overleftarrow{B})$.
3. $\sigma(\overrightarrow{B'}) = \sigma(\overleftarrow{A})$.
4. $\sigma(\overleftarrow{B'}) = \sigma(\overrightarrow{A})$.
5. $\sigma(\overleftarrow{C'}) = \sigma(\overleftarrow{B})$.
6. $\sigma(C' \to A') = [\sigma(\overleftarrow{B}) + 1] \cup [-\sigma(\overleftarrow{B}) + 1]$.
7. $\sigma(C' \to B') = \sigma(B \to A) + 1$.
8. $\sigma(B' \to A') = -\sigma(B \to A)$.
9. $\sigma(A' \to B') = \sigma(B \to A)$.

**Case 3.** *An Internal Block Creation move with parameter $b \geq 2$ is performed.*

Recall that this move can only be made if $a_1 > 2b_1$ and $B_b$ is a bottom block with left-most vertex $v_k$ where $k \leq a_1/2$ in $M$. Let $A'$ be the vertices in the newly created bottom block and let $B'$ be the vertices at the other end of top edges from $A'$ where we let $B$ be the set of vertices of $M$ corresponding to $B'$. Let $C$ and $D$ be the remaining vertices of $A_1$ on the left and right respectively of $B$ in $M$ and $C'$ and $D'$ be the corresponding vertices of $A'_1$ in $M'$. Note that no measures in or among the sets $B', C'$ or $D'$ have changed as a result of this move. Then the following represent all possible measures within the first top block of $M'$ as a result of the move.

1. $\sigma(\overrightarrow{A'}) = \sigma(\overleftarrow{B})$.
2. $\sigma(\overleftarrow{A'}) = -\sigma(\overleftarrow{B})$.
3. $\sigma(B' \to A') = [\sigma(\overleftarrow{B}) + 1] \cup [-\sigma(\overleftarrow{B}) + 1]$.
4. $\sigma(D' \to A') = \sigma(D \to B) - 1 = -\sigma(B \to C)$.
5. $\sigma(A' \to C') = \sigma(B \to C) + 1 = -\sigma(D \to B)$.

**Case 4.** *An Internal Rotation Expansion move is performed.*

Define sets in $M$ and $M'$ as in Figure 14. All measures within $C'$ and in between $B'$ and $C'$ in $M'$ are the same as those within $C$ and in between $B$ and $C$ in respectively in $M$. Also all measures within and among the sets $A', D'$ and $E'$ in $M'$ are the same as those within and among $A, B$ and $E$ in $M$. Thus, we consider all other pairs as represented in the following list.

1. $\sigma(\overrightarrow{B'}) = \sigma(\overleftarrow{B})$.

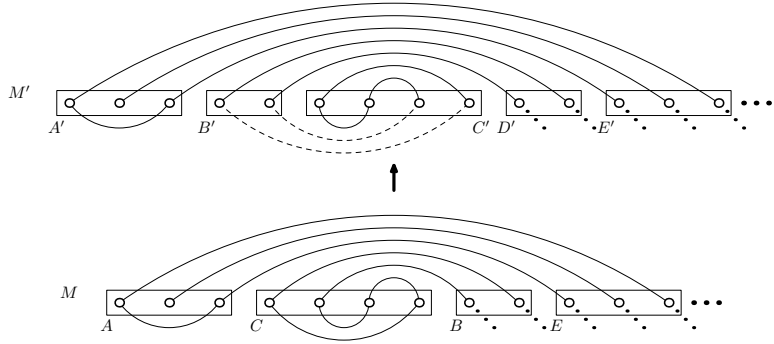

Figure 14: Internal Rotation Expansion.

2. $\sigma(\overleftarrow{B'}) = -\sigma(\overleftarrow{B})$.
3. $\sigma(D' \to B') = [\sigma(\overleftarrow{B}) + 1] \cup [-\sigma(\overleftarrow{B}) + 1]$.
4. $\sigma(D' \to C') = \sigma(B \to C) + 1$.
5. $\sigma(C' \to B') = -\sigma(B \to C)$.
6. $\sigma(E' \to B' \cup C') = \sigma(E \to B \cup C) + 1$.
7. $\sigma(B' \cup C' \to A') = \sigma(B \cup C \to A) - 1 = -\sigma(E \to B \cup C)$.

In all cases, the spectrum of $M'$ is symmetric and unbroken. Thus, the proof of the theorem is complete. $\square$

It appears that this proof technique might extend to establish the following conjecture.

**Conjecture 11.** *The dimensions are unimodal, that is, they are monotone increasing for values up to $0$ and monotone decreasing for values at least $1$.*

A proof of this conjecture using these techniques would likely require a deeper understanding of the moves. The complication lies in situations such as the one detailed in Item 7 of Case 2. Unimodality means that the measures are, in some sense, clustered near the center (at 0.5). However in this case, the measures contributed by the pairs between $B$ and $A$ may not be centered around 0.5. The addition of 1 to these measures may drive the set of measures further from the center. Thus, the addition of a balancing set which ensures the total is still symmetric may, in theory, not be enough to ensure such unimodality.

As a further strengthening of Theorem 9, we propose this conjecture.

**Conjecture 12.** *The dimensions are strictly unimodal, that is, they are strictly increasing for values up to 0 and strictly decreasing for values at least 1.*

*4.3. Index of a meander with 4 blocks*

The determination of whether a seaweed subalgebra is Frobenius usually relies on a combinatorial argument. This section recasts known formulas for computing the index. Theorem 15 provides a new addition to this class of formulas and is a strengthening of a result of Coll et al. from [2].

**Theorem 13 (Elashvili [7]).** *A meander of type $\frac{a|b}{c}$ has index $\gcd(a,b)-1$.*

In order to prove our next result, we state the following easy fact that follows immediately from the Euclidean Algorithm.

**Fact 14.** *For integers a and b,*

$$\gcd(a,b) = \gcd(a, a+b).$$

In the style of Theorem 13, we show the following result.

**Theorem 15.** *A meander of type $\frac{d}{a|b|c}$ or type $\frac{a|b}{c|d}$ has index $\gcd(a+b, b+c)-1$.*

Certainly a Flip Move has no effect on the index of the meander so we may assume any meander with four blocks is in one of these forms or its flip, whichever is more convenient at the time.

*Proof:* Let $M$ be a meander with 4 blocks. The proof is by induction on the number of simplified Winding Down moves in the signature of $M$. Let $n$ be the number of moves in the signature before either $C_0(c)$ or $B_0$ is performed (which reduce to a smaller number of blocks). The base of this induction, when $n = 0$, is provided by Cases 1 and 2. We consider cases based on the first simplified Winding Down move performed in the signature of $M$. Let $M'$ be the result of this first Winding Down move.

**Case 1.** *A Component Elimination move $C_0(c)$ is performed.*

In this case, $M'$ has only two blocks so must look like $\frac{a}{a}$. This means $M$ must have type $\frac{c|a}{c|a}$. Thus, the index of this meander is clearly $a+c-1$ while the formula says the index should be $\gcd(c+a, c+a) - 1 = a+c-1$.

**Case 2.** *A Block Elimination move $B_0$ is performed.*

This case has two subcases. The first is where $M = \frac{a|b}{c|d}$ and $a = 2c$. Then $M' = \frac{c|b}{d}$. By Theorem 13, the index of $M'$ is $\gcd(b,c) - 1$. Since $a = 2c$, using Fact 14, this value is the same as $\gcd(a+b, b+c) - 1$ and the Block Elimination move preserves the index, the index of $M$ is $\gcd(a+b, b+c) - 1$ as claimed. The second case is where $M = \frac{d}{a|b|c}$ so $d = 2a$ and $M' = \frac{a}{b|c}$. By Theorem 13, the index of $M'$ is $\gcd(b,c) - 1$. Since $d = a + b + c = 2a$, by Fact 14, this value is the same as $\gcd(a+b, b+c) - 1$ and the same argument yields the index of $M$.

**Case 3.** *A Rotation Contraction move $R_0$ is performed.*

If $M = \frac{a|b}{c|d}$, then $M' = \frac{c|b}{(2c-a)|d}$. By induction, the index of $M'$ is $\gcd(c+b, b+2c-a) - 1$. By Fact 14, $\gcd(c+b, b+2c-a) - 1 = \gcd(a+b, b+c) - 1$ and since the index is preserved under the Rotation Contraction move, the index of $M$ is as claimed. If $M = \frac{d}{a|b|c}$, an identical argument holds.

**Case 4.** *A Pure Contraction move $P_0$ is performed.*

If $M = \frac{a|b}{c|d}$, then $M' = \frac{(a-2c)|c|b}{d}$. By induction, the index of $M'$ is $\gcd(a-c, b+c) - 1$. By Fact 14, this equals $\gcd(a+b, b+c) - 1$ and since the index is preserved under the Pure Contraction move, the index of $M$ is as claimed. If $M = \frac{d}{a|b|c}$, an identical argument holds. $\square$

Note that Theorem 13 can be reproven very easily using exactly the same approach. A simple corollary of this result is the following theorem of Coll et al.

**Theorem 16 (Coll et al. [2]).** *A meander of type $\frac{d}{a|b|c}$ or type $\frac{a|b}{c|d}$ is Frobenius if and only if $\gcd(a+b, b+c) = 1$.*

*4.4. New infinite families of Frobenius meanders*

The known infinite families of Frobenius meanders are few. Frobenius meanders may have an arbitrarily large number of blocks, as in the classical examples Panyushev [16]. The following theorem shows that we can in fact have arbitrarily many blocks of arbitrarily large size. This can be easily extended to biparabolic meanders as seen in Corollary 18

**Theorem 17 (Coll et al. [3]).** *If $a$ is even and $\gcd(a,b) = 1$, then the meander of type $\frac{a|a|\cdots|a|b}{c}$ is Frobenius where $c$ is the sum of the terms on the top.*

The simplified Winding Down process can be used to obtain the following more general result.

**Corollary 18.** *If $a$ is even and $\gcd(a,b) = 1$, then the meander of type $\frac{a|a|\cdots|a|b}{c|a|a|\cdots|a}$ is Frobenius where $c = b + ka$ for some integer $k$.*

*Proof:* Let $a, b$ and $c$ be as given and let $\ell$ be the number of blocks of size $a$ on the bottom in this meander. By Theorem 17, the meander $M = \frac{a|a|\cdots|a|b}{(c+2\ell a)}$ is Frobenius (where there are $\ell$ extra blocks of size $a$ on the top of this meander). To $M$, we apply a simplified Flip $F_0$ followed by $\ell$ simplified Pure $P_0$ moves. The resulting meander is then flipped again using the $F_0$ move to obtain the meander $M' = \frac{a|a|\cdots|a|b}{c|a|a|\cdots|a}$. Since $M$ was Frobenius, $M'$ is as well. □

## 5. Conclusion

The classification of Frobenius Lie algebras appears to be a wild problem and the difficulty of classification may descend to seaweed Lie subalgebras of $\mathfrak{sl}(n)$. The wildness of the latter classification seems to present itself very quickly as the number of blocks increases. Indeed, following Theorems 13 and 15, one might expect that a Frobenius meander of type $a|b|c|d$ can be characterized by a relatively prime condition of the form

$$\gcd(\alpha_1 a + \alpha_2 b + \alpha_3 c + \alpha_4 d, \beta_1 a + \beta_2 b + \beta_3 c + \beta_4 d) = 1,$$

where the $\alpha_i$ and $\beta_j$ are integer *coefficients*. Substantial empirical evidence suggests that this is not so.

Exhaustive simulations have shown that there is no set of coefficients, all with absolute value at most ten, that can be used in such a relatively prime condition. All conditions with coefficients between negative ten and positive ten were checked against a set consisting of a large set of Frobenius meanders with five blocks. No condition survived. Since the addition of blocks seems to only complicate the situation, we are led to the following conjecture.

**Conjecture 19.** *There is no single relatively prime condition that suffices to classify Frobenius meanders with at least five blocks.*

Of course, if Conjecture 19 is true, it does not eliminate the possibility that some finite set of relatively prime conditions might classify meanders with five or more blocks that are Frobenius. However, (less than exhaustive) simulations of numerous sets of conditions were tested against millions of Frobenius meanders built up using the simplified Winding Up procedure outlined above. No set survived. We are led to the following stronger conjecture.

**Conjecture 20.** *There is no finite set of relatively prime conditions that suffices to classify Frobenius meanders with at least five blocks.*

In [15], Ooms shows by example that generally the eigenvalues of the adjoint of a principal element need not be integers. However, if the containing Lie algebra is semisimple and more generally algebraic, then Dergachev (unpublished) has observed that the eigenvalues are integers. We have shown that if $\mathfrak{g} = \mathfrak{sl}(n)$, then these eigenvalues form an unbroken sequence of integers where the dimensions of the associated eigenspaces are symmetric about $1/2$. More seems to be true.

**Conjecture 21.** *For $\mathfrak{g} = \mathfrak{sl}(n)$, the dimensions are strictly unimodal. In other words, no two consecutive terms in the sequence of dimensions are the same except for the dimensions of the eigenspaces of $0$ and $1$.*

An understanding of the conditions on $\mathfrak{g}$ (and $\mathfrak{f}$ with $\mathfrak{f} \subseteq \mathfrak{g}$) which ensure that the sequence of eigenvalues is an unbroken string of integers would be of interest. Subsequent modality questions naturally follow.

*Acknowledgement*

The authors would like to thank Murray Gerstenhaber & Anthony Giaquinto for numerous constructive conversations and suggestions regarding this work. We would also like to thank Anthony Giaquinto & Aaron Lauve for pointing out inconsistencies in a previous version of this work.**References**

[1] A. A. Belavin and V. G. Drinfel'd. Solutions of the classical Yang–Baxter equations for simple Lie algebras. *Funct. Anal. Appl.*, 16:159–180, 1982.25

[2] V. Coll, A. Giaquinto and C. Magnant. Meander graphs and Frobenius seaweed Lie algebras. *J. Gen. Lie Theory*, 5:Article ID G110103, 5 pages, 2011.

[3] V. Coll, C. Magnant and H. Wang. A new family of Frobenius parabolic Lie algebras. *Submitted*.

[4] V. Dergachev and A. Kirillov. Index of Lie algebras of seaweed type. *J. Lie Theory*, 10:331–343, 2000.

[5] M. Duflo. Personal communication to Gerstenhaber and Giaquinto.

[6] A. Dvorsky. Index of parabolic and seaweed subalgebras of $\mathfrak{so}_n$. *Linear Algebra and its App.*, 374:127–142, 2003.

[7] A. G. Elashvili. On the index of parabolic subalgebras of semisimple Lie algebras. preprint.

[8] A. G. Elashvili. Frobenius Lie algebras. *Funktsional. Anal. i Prilozhen.*, 16:94–95, 1982.

[9] M. Gerstenhaber and A. Giaquinto. The principal element of a Frobenius Lie algebra. *Letters Math. Phys.*, 88:333–341, 2009.

[10] M. Gerstenhaber and A. Giaquinto. Boundary solutions of the classical Yang-Baxter equation. *Letters Math. Phys.*, 40:337–353, 1997.

[11] M. Gerstenhaber and A. Giaquinto. Graphs, Frobenius functionals, and the classical Yang-Baxter equation. arXiv:0808.2423v1 [math.QA], August, 2008.

[12] A. Giaquinto, A. Lauve and J. Versnel. *Personal communication*.

[13] A. Joseph. On semi-invariants and index for biparabolic (seaweed) algebras, I. *J. of Algebra*, 305:487–515, 2006.

[14] S. M. Khoroshkin, I. I. Pop, M. E. Samsonov, A. A. Stolin, and V. N. Tolstoy, On some Lie bialgebra structures on polynomial algebras and their quantization. *Comm. Math. Phys.*, 282 (2008), 625662

[15] A. I. Ooms. On Lie algebras having a primitive universal enveloping algebra. *J. Algebra*, 32:488–500, 1974.

[16] D. Panyushev. Inductive formulas for the index of seaweed Lie algebras. *Mosc. Math. Journal*, 2:2001, 221-241.

[17] A. Stolin. On rational solution of Yang-Baxter equation for $\mathfrak{sl}(n)$. *Math Scand.*, 69 (1991), 5780.

[18] P. Tauvel and R. W. T. Yu. Sur l'indice de certaines algebres de Lie. *Ann. Inst. Fourier (Grenoble)*, 54:2004, 1793-1810.